\documentclass[10pt]{article}
\usepackage{amstext,amssymb,amsmath,amsbsy}

\usepackage{hyperref}
\usepackage{amscd}
\usepackage{amsfonts}
\usepackage{indentfirst}
\usepackage{verbatim}
\usepackage{graphicx}
\usepackage{color}

\usepackage{amsmath}
\usepackage{amsthm}
\usepackage{enumerate}
\usepackage{graphicx}
\usepackage[OT1]{fontenc}
\usepackage[latin1]{inputenc}
\usepackage[english]{babel}
\usepackage{amssymb}
\newtheorem{theorem}{Theorem}
\newtheorem{lemma}{Lemma}

\newtheorem{remark}{Remark}

\newtheorem{definition}{Definition}
\numberwithin{equation}{section}

\newcommand{\proofend}{\hfill $\Box$ }

\newcommand{\dsp}{\displaystyle}

\newcommand{\supp}{\operatorname{supp}}

\newcommand{\dive}{\operatorname{div}}

\newcommand{\eps}{\varepsilon}

\newcommand{\loc}{_{loc}}

\newcommand{\mR}{\mathbb{R}}

\newcommand{\mC}{\mathbb{C}}

\numberwithin{equation}{section}

\title{Localized and complete resonance in plasmonic structures}
\author{Hoai-Minh Nguyen \footnote{EPFL SB MATHAA CAMA, Station 8,  CH-1015 Lausanne, hoai-minh.nguyen@epfl.ch} \footnote{School of Mathematics, University of Minnesota, MN, 55455, hmnguyen@math.umn.edu} \footnote{The research is supported by NSF grant DMS-1201370 and by the Alfred P. Sloan Foundation.} \, and Loc Hoang Nguyen\footnote{EPFL SB MATHAA CAMA, Station 8,  CH-1015 Lausanne, loc.nguyen@epfl.ch}}

\begin{document}
\maketitle

\begin{abstract}
This paper studies a possible  connection between the way the  time averaged electromagnetic power dissipated into heat blows up and the anomalous localized resonance in plasmonic structures. We show that there is a setting in which the localized resonance takes place whenever the resonance does and moreover, the power is always bounded and might go to $0$.  We also provide another setting in which the resonance is complete and the power goes to infinity whenever resonance occurs; as a consequence of this fact there is no localized resonance.   This work is motivated from recent works on cloaking via anomalous localized resonance.
\end{abstract}

\section{Introduction and statement of the main results}
Negative index materials (NIMs) were first investigated theoretically by Veselago in \cite{Veselago} and were innovated by Nicorovici et al. \cite{NicoroviciMcPhedranMilton94} in the electrical impedance setting and by Pendry  \cite{PendryNegative} in the electromagnetic setting. The existence of such materials was confirmed by Shelby, Smith, and Schultz in \cite{ShelbySmithSchultz}. An interesting (and surprising) property on NIMs is the anomalous localized resonance discovered by Nicorovici et al.  in \cite{NicoroviciMcPhedranMilton94} for core-shell plasmonic structures in two dimensions in which a circular shell has permitivity $-1+i\delta$ while the core and the matrix, 
the complement of the core-shell structure, have permitivity $1$. Here $\delta$ describes the loss of the material (more precisely, the loss of the negative index material part). 
A key figure of the phenomenon is the  localized resonance of the field, i.e., 
the field blows up in some regions and remains bounded in some others as $\delta \to 0$. This is  partially due to the change sign of the coefficient in the equation and therefore the ellipticity is lost as $\delta \to 0$;  the loss of ellipticity  is not sufficient to ensure such a property as discussed later in this paper. Following \cite{MiltonNicorovici}, the localized resonance  is anomalous because  the boundary of the resonant regions varies with the position of the source, and their boundary does not coincide with any discontinuity in moduli. 

\medskip
An attractive application related to the anomalous localized resonance is cloaking. This was recognized by Milton and Nicorovici in \cite{MiltonNicorovici} and investigated in \cite{AmmariCiraoloKangLee, AmmariCiraoloKangLeeMilton2, AmmariCiraoloKangLeeMilton1, BouchitteSchweizer10, KohnLu, Milton-folded} and the references therein.  Let us discuss two results related to cloaking via anomalous localized resonance  obtained so far for non radial core shell structures  in \cite{AmmariCiraoloKangLee, KohnLu},  in which the authors deal with  the two dimensional quasistatic regime.  In \cite{AmmariCiraoloKangLee}, the authors provide a necessary and sufficient condition on the  source for which the  time averaged electromagnetic power dissipated into heat blows up as the loss  goes to zero using the spectral method. Their characterization is based on the  detailed information on the spectral properties of a Neumann-Poincar\'e type operator. This information is difficult to come by in general. In \cite{KohnLu}, using the variational approach, the authors show that the power goes to infinity if the location of the source is in a finite range w.r.t. the shell for a class of sources. The core is not assumed to be radial  but the matrix is in \cite{KohnLu}.  The boundedness of the fields in some regions for these structures is not discussed in \cite{AmmariCiraoloKangLee, KohnLu}  except in the radial case showed in \cite{AmmariCiraoloKangLee} (see also \cite{MiltonNicorovici, NicoroviciMcPhedranMilton94}).  It is of interest to understand if there is a possible connection between  the power and the localized resonance in general. 

\medskip
In this paper, we present two settings  in which there is no connection between the blow up of  the power and the localized resonance. To this end,  the following two problems are considered. 

\medskip

\noindent{\bf Problem 1:} The behaviour of $u_\delta \in H^1(B_R)$ ($R>1$) the unique solution to 
\begin{equation}\label{eq-P2}
\left\{\begin{array}{cl}
\dive(\eps_\delta \nabla u_\delta) = 0 & \mbox{ in } B_R,  \\[6pt]
u_\delta = g & \mbox{ on } \partial B_R, 
\end{array} \right.
\end{equation}
where $g \in H^{1/2}(\partial B_R)$ and the way the power, which will be defined in \eqref{def-E}, explodes as $\delta \to 0_+$. 

\medskip 

Here and in what follows $B_r$ denotes the ball centred at the origin of radius $r$ for $r > 0$.

\medskip

\noindent {\bf Problem 2:} The behaviour of $u_\delta \in W^1(\mR^2)$ (see \eqref{def-W1} for the notation) the unique solution converging to 0 as $|x| \to \infty$ to 
\begin{equation}\label{eq-P1}
\dive(\eps_\delta \nabla u_\delta) = f \mbox{ in } \mR^2, 
\end{equation}
and the way the power, defined in \eqref{def-E},  explodes. Here
$f$ is  in $L^2(\mR^2)$ with compact support in $\mR^2 \setminus B_1$ and  satisfies the  compatible condition
\begin{equation}\label{cond1}
\int_{\mR^2} f = 0. 
\end{equation}

For $0 \le  \delta < 1$,  $\eps_\delta$ is defined by
\begin{equation}\label{def-eps}
\eps_\delta : = \left\{ \begin{array}{cl} (F^{-1})_{*}I & \mbox{ if } |x| > 1\\[6pt]
-1 + i \delta & \mbox{ if } |x| < 1, 
\end{array}\right. \quad \mbox{ for } d= 2, \, 3, 
\end{equation}
where $F: \mR^{d} \setminus B_1 \to \bar B_1$ is the Kelvin transform w.r.t. $\partial B_{1}$, i.e., $F(x) = x/ |x|^{2}$. 

\medskip
Here and in what follows, we use the following standard notation
\begin{equation}\label{defTO}
T_*a(y) = \frac{D T  (x)  a(x) D T ^{T}(x)}{J(x)} \quad \mbox{ and } \quad T_*f (y) = \frac{f (x)}{J(x)}. 
\end{equation}
where $x =T^{-1}(y)$ and $J(x) = |\det D T(x)|$, for $f \in L^2(D_1)$, $a \in [L^\infty(D_1)]^{d \times d}$, $f \in L^2(D_1)$, and $T$ a diffeomorphism from $D_1$  onto $D_2$.

\medskip 
It is easy to verify that, as noted in \cite{Ng-Complementary},  
 \begin{equation*}
\eps_\delta : = \left\{ \begin{array}{cl} 1 & \mbox{ if } |x| > 1\\[6pt]
-1 + i \delta & \mbox{ if } |x| < 1
\end{array}\right. \quad \mbox{ for } d=2. 
\end{equation*}  

The media considered in Problems 1 and 2 where $\eps_\delta$ is given in \eqref{def-eps} have the complementary property (see \cite{Ng-Complementary} for the definition and a discussion on various results related to these media in a general core shell structure). 
The setting studied in \cite{KohnLu} also inherits this property since the matrix is radial while the setting in \cite{AmmariCiraoloKangLee} is not in general.  As seen later, this property is not enough to ensure a connection between the blow up of the power and the localized resonance. 

\medskip
In Problems~1 and 2, $\delta$ is the loss of the media (more precisely the loss of the negative index material in $B_1$) and the time averaged power dissipated into heat is given by (see, e.g.,  \cite{AmmariCiraoloKangLee, KohnLu})
\begin{equation}\label{def-E}
E_\delta (u_\delta) = \delta \int_{B_1} |\nabla u_\delta|^2 \, dx. 
\end{equation}
From the definition of $u_\delta$, one can derive that 
\begin{equation*}
 \int_{B_1} |\nabla u_\delta|^2 \, dx \ge \left\{ \begin{array}{cl} \dsp C_1  \int_{B_R} \dsp |\nabla u_\delta|^2 - C_2  \|g\|_{H^{1/2}(\partial B_R)}^2 & \mbox{ in  Problem 1}, \\[6pt]
\dsp  C_1 \int_{\mR^2} |\nabla u_\delta|^2 - C_2 \| f\|_{L^2}^2 & \mbox{ in Problem 2}, 
\end{array}\right. 
\end{equation*}
for some positive constants $C_1$, $C_2$ independent of $\delta$, $f$, and $g$.

\medskip
The main results of the paper are Theorems~\ref{thm2} and \ref{thm1} below. Concerning Problem 1, we have.

\begin{theorem} \label{thm2} Let  $d=2, 3$, and $g \in H^{1/2}(\partial B_R)$ and $u_\delta \in H^1(B_R)$ be the unique solution to \eqref{eq-P2}. Then 
\begin{enumerate}
\item Case 1: $g$ is compatible to \eqref{eq-P2} (see Definition~\ref{def1}). Then $\big(\|u_\delta\|_{H^1(B_R)} \big)$ remains bounded. Moreover, $u_\delta \to u_0$ weakly in $H^1(B_R)$ as $\delta \to 0$ where $u_0 \in H^1(B_R)$ is the unique solution to  
\begin{equation}\label{state-1-P2}
\left\{\begin{array}{cl}
\dive(\eps_0 \nabla u_0) = 0 & \mbox{ in } B_R, \\[6pt]
u_0 = g & \mbox{ on } \partial B_R.  
\end{array} \right.
\end{equation}

\item Case 2: $g$ is not compatible to  \eqref{eq-P2}. Then 
\begin{equation}\label{state-2-P2}
\lim_{\delta \to 0} \|u_\delta\|_{H^1(B_R)} = + \infty; 
\end{equation}
however,  
\begin{equation}\label{state-3-P2}
u_\delta \to v \mbox{ weakly in } H^1(B_{1/R}), 
\end{equation}
where $v \in H^{1}(B_{1/R})$ is the unique solution to 
\begin{equation}\label{def-v}
	\left\{
		\begin{array}{cl}
			\Delta v = 0 &\mbox{in } B_{1/R},\\[6pt]
			v(x) = h(x): = g(x/|x|^{2})  &\mbox{on } \partial B_{1/R}. 
		\end{array}
	\right.
\end{equation} 
Moreover, for all $g \in H^{1/2}(\partial B_R)$, 
\begin{equation}\label{state-4-P2}
\limsup_{\delta \to 0} \delta \int_{B_R} |\nabla u_\delta|^2 dx < + \infty, 
\end{equation}
and for any $0 < \alpha < 1/2$, there exists $g \in H^{1/2}(\partial B_R)$ such that 
\begin{equation}\label{state-5-P2}
0 <  \liminf_{\delta \to 0} \delta^{2 \alpha} \int_{B_R} |\nabla u_\delta|^2 dx \le \limsup_{\delta \to 0} \delta^{2 \alpha} \int_{B_R} |\nabla u_\delta|^2 dx < + \infty. 
\end{equation}
\end{enumerate}
\end{theorem}

\begin{remark}\label{rem-P1}
Concerning~\eqref{eq-P2}, whenever resonance takes place \footnote{In \cite{AmmariCiraoloKangLee} and \cite{KohnLu}, the authors introduced the definition of resonance. Following them, a system is resonant if and only if the power blows up as $\delta \to 0$. }, it is localized  in the sense that the field blows up in some region and remains bounded in some others; moreover, the power remains bounded and might converge to 0 as $\delta \to 0$ \footnote{Graeme Milton recently informed us that some examples on anomalous localized resonance (for dipole sources) without the blow up of the power are given in \cite{MiltonNicoroviciMcPhedranPodolskiy}. We thank him for pointing this out. We note here that the setting in this paper is different from the one in \cite{MiltonNicoroviciMcPhedranPodolskiy} where the negative index material part is  in a shell not in a ball; the anomalous localized resonance and boundedness of the power in the setting in \cite{MiltonNicoroviciMcPhedranPodolskiy} depend on the location of the source.}. 
\end{remark}

In the statement of Theorem~\ref{thm2}, we use the following definition. 

\begin{definition}\label{def1}
	Let $g \in H^{1/2}(\partial B_R)$. Then $g$ is said to be {\bf compatible} to \eqref{eq-P2} if and only if there exists a solution $w \in H^1(B_{1} \setminus  B_{1/R})$ to the Cauchy problem
\begin{equation}
	\left\{
		\begin{array}{cl}
			\Delta w = 0 &\mbox{in } B_{1} \setminus B_{1/R},\\[6pt]
			w = v \mbox { and } \partial_{\nu} w = \partial_{\nu} v  &\mbox{on } \partial B_{1/R},  
		\end{array}
	\right.
\label{eqn compatible in BR}
\end{equation} 
where $v$ is the function defined in \eqref{def-v}. Otherwise, $g$ is not compatible. 
\end{definition}

\begin{remark} Figure \ref{fig 1} in Section \ref{Sec Numer} provides a numerical simulation illustrating Theorem~\ref{thm2}. 
\end{remark}


Concerning Problem~2, we have. 

\begin{theorem}\label{thm1} Let $f \in L^2(\mR^2)$ be such that $\supp f \subset \subset  \mR^d \setminus B_1$ and  \eqref{cond1} holds and let $u_\delta \in W^1(\mR^2)$ be the unique solution converging to 0 as $|x| \to \infty$ to \eqref{eq-P1}. Then
\begin{enumerate}
\item Case 1: $f$ is compatible to  \eqref{eq-P1} (see Definition~\ref{def2}). Then 
\begin{equation}\label{state-1-P1}
u_\delta = U: = \left\{\begin{array}{cl} w \circ F - w(0) &  \mbox{ in } \mR^{2} \setminus B_{1}, \\[6pt]
- w(0) & \mbox{ in } B_{1}. 
\end{array} \right.
\end{equation}  
Here $w$ will be defined in \eqref{def-w-P1}. 
 
\item Case 2: $f$ is not compatible to  \eqref{eq-P1}. Then 
\begin{equation}\label{state-2-P1}
0 < \liminf_{\delta \to 0} \delta^2 \int_{O} |\nabla u_\delta|^2 \, dx  \leq \limsup_{\delta \to 0} \delta^2 \int_{O} |\nabla u_\delta|^2 \, dx < + \infty, 
\end{equation}
for any open subset $O$ of $\mR^2$. 
\end{enumerate}
\end{theorem}

\begin{remark}
	Inequalities \eqref{state-2-P1} implies that the field blows up in any open subset of $\mathbb{R}^2$ at the same rate \footnote{Graeme Milton recently informed us that for a single dipole source outside $B_1$,  the resonance is not localized.}.
\end{remark}

\begin{remark} Theorem~\ref{thm1} also holds for $d=3$  (see the proof of Theorem~\ref{thm1} and Remark~\ref{rem-proof-comments}, which is about representations in $B_1$). However, in this case, the existence of $u_\delta$ belongs to some Sobolev spaces with weight since $(F^{-1})_* I$ is not bounded from below by a positive constant at infinity due to the fact $d=3$. We do not treat this case in this paper to keep the presentation simple. 
\end{remark}

For $U$ a smooth open region of $\mR^2$ with a bounded complement (this includes $U = \mR^2$), we use the following standard notation:
\begin{equation}\label{def-W1}
	W^1(U) = \left\{u \in L^2_{\loc}(U);  \;  \nabla u \in \big[L^2(U)\big]^2  \mbox{ and } \frac{u}{|x|\log(2 + |x|)} \in L^2(U)\right\}.
\end{equation}

Part of Theorem~\ref{thm1} was considered in \cite{KohnLu}. More precisely, in \cite{KohnLu}, the authors showed that $E_\delta (u_\delta) \to \infty$ for $f$ with $\supp f \subset \partial B_r$ for $r > 1$ \footnote{In fact, such an $f$ is not in $L^2(\mR^2)$, however our analysis is  also valid  for this case. Our presentation is restricted for $f \in L^2$ so that the definition of $\big(F^{-1}\big)_*f$ makes sense without introducing further notations.}. In this paper, we make one step further. We show that when resonance occurs, it is {\bf complete} in the sense that \eqref{state-2-P1} holds;  there is no localized resonance here. Otherwise, the field remains bounded.  In fact it is independent of $\delta$ by \eqref{state-1-P1}.

\medskip 
In the statement of Theorem~\ref{thm1}, we use the following definition. 

\begin{definition} \label{def2}
	Let $f \in L^2(\mR^2)$ with $\supp f \subset  \mR^2 \setminus B_1$. Then $f$ is said to be \textbf{compatible} to  \eqref{eq-P1} if and only if there exists a solution $w \in H^1(B_1)$ to the Cauchy problem
	\begin{equation}\label{def-w-P1}
	\left\{
		\begin{array}{cl}
			\Delta w = F_{*}f & \mbox{ in } B_{1},\\[6pt]
			\partial_{\nu} w  = w  = 0 & \mbox{ on } \partial B_1.  
		\end{array}
	\right.
\end{equation}
Otherwise, $f$ is not compatible. 
\end{definition}

\begin{remark} Figure \ref{fig 2} in Section \ref{Sec Numer} provides a numerical simulation illustrating Theorem~\ref{thm1}. 
\end{remark}




From Theorems~\ref{thm2} and \ref{thm1}, we conclude that in the settings considered in this paper, there is no connection between the unboundedness of the  power and the localized resonance. Though the settings in Problems 1 and 2 are very similar, the essence of the resonance  are very different.  A connection between these phenomena would be linked  not only to the location of the source but also to the geometry of the problem, i.e., the definition of $\eps_\delta$.  Using the concept of (reflecting) complementary media introduced in \cite{Ng-Complementary}, one can extend the results this paper in a more general setting.

\medskip 
The definitions of compatibility conditions have roots from \cite{Ng-Complementary}.  The analysis for the compatible cases is inspired from there. The analysis in the incompatible case is guided from the compatible one. One of the main observations in this paper is the localized resonant phenomena in \eqref{state-3-P2} (one has localized resonance by \eqref{state-2-P2}). The localized resonance is also discussed in the context of  superlensing and cloaking using  complementary media  in \cite{Ng-superlensing, Ng-Negative-Cloaking} where the removing of localized singularity technique was introduced by the first author to deal with localized resonance in non radial settings.  In recent work \cite{Ng-CALR}, the first author introduces the concept of doubly complementary media for a general shell-core structure and shows that cloaking via anomalous localized resonance takes place if and only if the power blows up. To this end, he  introduces and develops the technique of separation of variables for a general structure. 
  
\medskip
The paper is organized as follows. In Sections~\ref{sect-thm2} and \ref{sect-thm1}, we prove Theorems~\ref{thm2} and \ref{thm1} respectively. 
In Section \ref{Sec Numer}, we provide numerical simulations illustrating  these results.

\section{Proof of Theorem~\ref{thm2}} \label{sect-thm2} 

\subsection{Preliminaries}
 In this section, we present two elementary lemmas which are very useful for the proof of Theorem~\ref{thm2}. The first one (Lemma~\ref{lem1}) is on the change of variables for the Kelvin transform. Lemma~\ref{lem1} is a special case of \cite[Lemma 4]{Ng-Complementary} which deals with general reflections. 
 
 \begin{lemma}\label{lem1} Let $d =2, \, 3$,  $0 < R_1 < R_2 < R_3$ with $R_3 = R_2^2/ R_1$, $f \in L^2(B_{R_2} \setminus B_{R_1})$,  $a \in [L^\infty(B_{R_2 \setminus R_1})]^{d \times d}$ be a uniformly elliptic  matrix - valued  function, and $K:B_{R_2} \setminus \bar B_{R_1} \to B_{R_3} \setminus \bar B_{R_2}$ be the Kelvin transform w.r.t $\partial B_{R_2}$,  i.e., 
\begin{equation*}
K(x) = R_2^2 x/ |x|^2. 
\end{equation*}
For $v \in H^1(B_{R_2} \setminus B_{R_1})$, define $w = v \circ K^{-1}$. Then 
\begin{equation*}
\dive (a \nabla v) = f \mbox{ in }  B_{R_2} \setminus B_{R_1}
\end{equation*}
if and only if
\begin{equation*}
\dive (K_*a \nabla w)  = K_*f \mbox{ in } B_{R_3} \setminus B_{R_2}.  
\end{equation*}
Moreover, 
\begin{equation*}
w = v \quad \mbox{ and } \quad K_*a \nabla w \cdot \nu = - a \nabla v \cdot \nu  \mbox{ on } \partial B_{R_2}.    
\end{equation*}
\end{lemma}
 
The second lemma is on an estimate related to solutions to \eqref{eq-P2}.  
\begin{lemma}\label{lem2}
Let $d=2, \, 3$, $f \in H^{-1}( B_R)$,  and  let $U_\delta \in H^1_0(B_R)$ be  the unique solution to 
\begin{equation*}
\dive (\eps_\delta \nabla U_\delta) = f \mbox{ in } B_R. 
\end{equation*}
We have
\begin{equation*}
\| U_\delta\|_{H^1(B_R)} \le \frac{C}{\delta} \| f\|_{H^{-1}}
\end{equation*}
for some positive constant $C$ independent of $f$ and $\delta$. 
\end{lemma}

\noindent {\bf Proof.} Lemma~\ref{lem2} follows from  Lax-Milgram's theorem. The details are left to the reader. \proofend

\subsection{Proof of Theorem~\ref{thm2}} \label{section Proof localized resonance}

The proof is divided into 6 steps. 

\medskip

\noindent \underline{Step 1:} We prove that  if
there exists a solution $u \in H^1(B_R)$ to 
\begin{equation}\label{eq-limit}\left\{\begin{array}{cl}
\dive(\eps_0 \nabla u) = 0 & \mbox{ in } B_R, \\[6pt]
u = g & \mbox{ on } \partial B_R,   
\end{array} \right.
\end{equation}
then $g$ is compatible. Moreover, the solution to \eqref{eq-limit} is unique in $H^1(B_R)$.  

\medskip
In fact,  define $V$ in $B_{1} \setminus B_{1/R}$ by
\begin{equation*}
V = u \circ F^{-1}.  
\end{equation*}
We have, by Lemma~\ref{lem1}, 
\begin{equation*}V = u \Big|_{\rm ext} \mbox{ and } \partial_{r} V =  \partial_{r} u \Big|_{\rm ext} \mbox{ on } \partial B_{1}. 
\end{equation*}
Set 
\begin{equation*}
W = u - V \mbox{ in } B_{1} \setminus B_{1/R}.
\end{equation*}
By Lemma~\ref{lem1},  $W \in H^{1}(B_{1} \setminus B_{1/R})$ is a solution to the Cauchy problem
\begin{equation*}
	\left\{
		\begin{array}{cl}
			\Delta W = 0 & \mbox{ in } B_{1} \setminus B_{1/R},\\[6pt]
			\partial_{\nu} W  = W  = 0 & \mbox{ on } \partial B_1.  
		\end{array}
	\right.
\end{equation*}
By the unique continuation principle, $W=0$. This implies 
\begin{equation*}
u = V= h \mbox{ on } \partial B_{1/R}. 
\end{equation*}
Therefore, $u = v$ in $B_{1/R}$ where $v$ is defined in \eqref{def-v}. 
It follows that  $u$ satisfies \eqref{eqn compatible in BR} and $g$ is compatible. The uniqueness in $H^1(B_R)$ of \eqref{eq-limit}  is also clear from the analysis. 

\medskip
\noindent \underline{Step 2:} We prove that if $g$ is compatible then $u_\delta \to u$ weakly in $H^1(B_R)$ where
\begin{equation*}
u=  \left\{\begin{array}{cl} 
v & \mbox{ in } B_{1/R}, \\[6pt]
w & \mbox{ in } B_1 \setminus B_{1/R}, \\[6pt]
 w \circ F & \mbox{ in } B_R \setminus B_1, 
\end{array} \right.
\end{equation*}
where $w$ is given in \eqref{eqn compatible in BR}. 

\medskip
It is clear that $u \in H^1(B_R)$ is a solution to \eqref{eq-limit}. The uniqueness of $u$ follows from Step 1. 
 Define
\begin{equation*}
U_\delta = u_\delta - u \quad \mbox{ in } B_R. 
\end{equation*}
Then $U_\delta \in H^1_0(B_R)$ is the unique solution to
\begin{equation*}
\dive (\eps_\delta \nabla U_\delta) = \dive \big( (\eps_0 - \eps_\delta) \nabla u \big) \mbox{ in } B_R. 
\end{equation*}
This implies, by Lemma~\ref{lem2},  
\begin{equation*}
\| U_\delta \|_{H^1(B_R)} \le C \|\nabla u \|_{L^2(B_R)}. 
\end{equation*}
Since $u_\delta = U_\delta + u$, $(u_\delta)$ is bounded in $H^1(B_R)$. W.l.o.g. one may assume that $u_\delta$ converges weakly in $H^1(B_R)$ to a solution to \eqref{eq-limit}. Since \eqref{eq-limit} is uniquely solvable in $H^1(B_R)$, the conclusion follows.  

\medskip
\noindent \underline{Step 3:} We prove that if $\liminf_{\delta \to 0} \| \nabla u_{\delta} \|_{L^{2}(\mR^{2})} < + \infty$ then $g$ is compatible. 

\medskip 
Since $\liminf_{\delta \to 0} \| \nabla u_{\delta} \|_{L^{2}(\mR^{2})} < + \infty$, there exists a solution $u \in H^1(B_R)$ to \eqref{eq-limit}. 
The conclusion now is a consequence of Step 1. 

\bigskip 

After Steps 1, 2, and 3, the first statement of Theorem~\ref{thm2} and \eqref{state-2-P2} are established. We next prove \eqref{state-3-P2}, \eqref{state-4-P2}, and \eqref{state-5-P2}. 
We will only consider the two dimensional case. The proof in three dimensions follows similarly (see Remark~\ref{rem-proof-comments}). In what follows, we assume that $d=2$. 

\bigskip
\noindent \underline{Step 4:} Proof of \eqref{state-3-P2}.

\medskip
Set 
\begin{equation*}
v_\delta = u_\delta \circ F^{-1} \mbox{ in } B_{1} \setminus B_{1/R}. 
\end{equation*}
Then $v_\delta \in H^1 (B_{1} \setminus B_{1/R})$ and 
\begin{equation*}
\Delta v_\delta = 0 \mbox{ in } B_{1} \setminus B_{1/R}. 
\end{equation*}
One can represent $v_\delta$ as follows 
\begin{equation}\label{form1}
v_\delta = a_0 + b_0 \log r + \sum_{n = 1}^{\infty} \sum_{\pm}(a_{n, \pm} r^n + b_{n, \pm} r^{-n}) e^{ \pm i n \theta} \mbox{ in } B_{1} \setminus B_{1/R}, 
\end{equation}
for  $a_0, \, b_0, a_{n, \pm}, b_{n, \pm} \in \mC$ ($n \ge 1$). 
Similarly,  one can represent $u_\delta$ by 
\begin{equation}\label{form2}
	u_{\delta} =  c_0 + \sum_{n = 1}^{\infty} \sum_{\pm} c_{n, \pm} r^n e^{\pm in \theta} \mbox{ in } B_{1}, 
\end{equation} 
for $c_0, \, c_{n, \pm} \in  \mC$ ($n \ge 1$). 
Using the transmission conditions on $\partial B_1$, we have
\begin{equation}\label{form3}
v_\delta = u_\delta \Big|_{\rm int} \quad \mbox{ and } \quad \partial_{\nu} v_\delta = (1 - i \delta) \partial_\nu u_\delta \Big|_{\rm int} \quad \mbox{ on } \partial B_1. 
\end{equation}
A combination of \eqref{form1}, \eqref{form2}, and \eqref{form3} yields 
\begin{equation*}
	\left\{ \begin{array}{c}
		a_{n, \pm} + b_{n, \pm} = c_{n, \pm}\\[6pt]
		a_{n, \pm} - b_{n, \pm} = (1 - i\delta)c_{n, \pm},
	\end{array}
	\right. \quad \mbox{ for  } n \ge 1,
\end{equation*}
and 
\begin{equation}\label{2.5}
	\left\{ \begin{array}{c}
		a_0  = c_0\\[6pt]
	           b_0 = 0.
	\end{array}
	\right.  
\end{equation}
This implies 
\begin{equation}\label{form1.2}
	\left\{ \begin{array}{l}
		\dsp a_{n, \pm} = (2 - i \delta)c_{n, \pm}/ 2 \\[6pt]
		\dsp b_{n, \pm} = i \delta c_{n, \pm} / 2
	\end{array}
	\right.  \mbox{ for } n \ge 1. \quad 
\end{equation}
From the definition of $v_\delta$, it is clear that 
\begin{equation}\label{form1.3}
v_\delta = h = h_0  + \sum_{n = 1}^\infty \sum_{\pm} h_{n, \pm}  e^{\pm i n \theta} \mbox{ on } \partial B_{1/R}, 
\end{equation}
for some $h_0, \, h_{n, \pm} \in \mC$ ($n \ge 1$). Since $v_\delta = h$ on $\partial B_{1/R}$, it follows from   \eqref{form1}, \eqref{2.5}, \eqref{form1.2}, and \eqref{form1.3} that
\begin{equation}\label{comp1}
\frac{1}{2} \left[ (2 - i \delta) R^{-n} + i \delta R^{n})\right]c_{n, \pm}  = h_{n, \pm} \mbox{ for } n \ge 1
\end{equation}
and 
\begin{equation}\label{comp2}
 c_0 = h_0. 
\end{equation}

We claim that 
\begin{equation}\label{claim}
u_\delta - v_\delta \to 0 \mbox{ weakly in } H^{1/2}(\partial B_{1/R}).
\end{equation}
In fact,  by \eqref{2.5} and \eqref{form1.2},  we have
\begin{equation}\label{comp3}
u_\delta - v_\delta = \sum_{n = 1}^\infty \sum_{\pm} \frac{1}{2} i \delta c_{n, \pm} (R^{-n} - R^{n}) e^{\pm i n \theta} \quad \mbox{ on } \partial B_{1/R}. 
\end{equation}
We derive from \eqref{comp1} and \eqref{comp3} that 
\begin{equation}\label{comp4}
u_\delta - v_\delta = \sum_{n =1}^\infty \sum_{\pm} \frac{ i \delta (R^{-n} - R^n) }{\big[ 2 R^{-n}  - i \delta (R^{-n} - R^n)\big]}  h_{n, \pm} e^{\pm i n \theta} \quad \mbox{ on } \partial B_{1/R}. 
\end{equation}
Claim \eqref{claim} follows since 
\begin{equation*}
\lim_{\delta \to 0} \frac{ i \delta (R^{-n} - R^n) }{\big[ 2 R^{-n}  - i \delta (R^{-n} - R^n)\big]} = 0 \mbox{ for all } n \ge 1
\end{equation*}
and 
\begin{equation*}
\left| \frac{ i \delta (R^{-n} - R^n) }{\big[ 2 R^{-n}  - i \delta (R^{-n} - R^n)\big]} \right| \le 1 \mbox{ for all } n \ge 1.  
\end{equation*}

The conclusion of Step 4 is now a consequence of Claim \eqref{claim} and the fact that $\Delta (u_{\delta} - v) = 0$ in $B_{1/R}$. 

\medskip 
\noindent \underline{Step 5:} Proof of \eqref{state-4-P2}: 

\medskip 
Since $\Delta u_\delta = 0$ in $B_R \setminus \partial B_1$ and $u_\delta = g$ on $\partial B_R$, it suffices to prove that
\begin{equation*}
\limsup_{\delta \to 0}  \delta \| u_\delta \|_{H^{1/2}(\partial B_1)}^2 \le C \| h \|_{H^{1/2}(\partial B_{1/R})}^2.  
\end{equation*}
In this proof, $C$ denotes a positive constant independent of $\delta$ and $g$.  From \eqref{form2}, 	 \eqref{comp1}, and \eqref{comp2}, we have
\begin{equation*}
	C \|u_{\delta}\|_{H^{1/2}(\partial B_1)}^2 \le   |h_0|^2 +  \sum_{n = 1}^{\infty}  \sum_{\pm} \frac{ n |h_n|^2}{4R^{-2n}  + \delta^2 (R^n - R^{-n})^2 }.
\end{equation*}
We derive that
\begin{equation}\label{tt1}
	C \|u_{\delta}\|_{H^{1/2}(\partial B_1)}^2  \le  \sup_{n \ge 0}  \frac{ 1}{4R^{-2n} + \delta^2 (R^n - R^{-n})^2 }   \Big( |h_0|^2 + \sum_{n = 1}^{\infty} \sum_{\pm} n |h_{n, \pm}|^2 \Big) . 
\end{equation}
Since $R>1$, it follows that  
\begin{equation}\label{tt2}
R^{-2n} + \delta^2 (R^n - R^{-n})^2 \ge C( R^{-2n} + \delta^2 R^{2n} ) \ge 2 C \delta \quad \forall \, n \ge 1. 
\end{equation}
A combination of  \eqref{tt1} and \eqref{tt2} yields 
\begin{equation*}
C \delta \|u_{\delta}\|_{H^{1/2}(\partial B_1)}^2 \le \| h \|_{H^{1/2}(\partial B_{1/R})}^2;  
\end{equation*}
hence \eqref{state-4-P2} follows. 

\medskip
\noindent \underline{Step 6:} Proof of \eqref{state-5-P2}. 

Since $\Delta u_\delta = 0$ in $B_R \setminus \partial B_1$ and $u_\delta = g$ on $\partial B_R$, it suffices to find $h \in H^{1/2}(\partial B_{1/R})$ such that  
\begin{equation}\label{tt34}
0 <  \liminf_{\delta \to 0} \delta^{2 \alpha} \| u_\delta\|_{H^{1/2}(\partial B_1)}^2 \le \limsup_{\delta \to 0} \delta^{2 \alpha}  \| u_\delta\|_{H^{1/2}(\partial B_1)}^2 < + \infty.
\end{equation}
Recall that $h(x) = g(x/ |x|^2)$. Let  $n_{\delta} = [\frac{1}{2}|\ln \delta/ \ln R|]$ be the smallest integer that is greater than or equal to $\frac{1}{2}|\ln \delta/ \ln R|$ ($R^{-2 n_\delta} \sim \delta$). We have 
\begin{align*}
	\|u_{\delta}\|_{H^{1/2}(\partial B_1)}^2 \sim&  |h_0|^2 + \sum_{n = 1}^{\infty} \sum_{\pm}  \frac{ n|h_{n, \pm}|^2}{4R^{-2n} + \delta^2 (R^n - R^{-n})^2 }\\[6pt]
	\sim & |h_0|^2 + \sum_{n = 1}^{n_{\delta}} \sum_{\pm} \frac{ n|h_{n, \pm}|^2}{R^{-2n}  } + \sum_{n = n_{\delta} + 1}^{\infty}  \sum_{\pm} \frac{ n |h_{n, \pm}|^2}{\delta^2 R^{2n}}. 
\end{align*}
Set
\begin{equation*}
0 < \gamma = 1 - 2\alpha  < 1
\end{equation*}
and choose 
\begin{equation*}
h_0 = 0 \quad \mbox{ and } \quad h_{n, \pm} = \frac{R^{-n \gamma}}{\sqrt{n}} \mbox{ for } n \ge 1. 
\end{equation*}
It follows that, since $\gamma < 1$ and $R>1$, 
\begin{equation}\label{tt3}
|h_0|^2 + 	\sum_{n = 1}^{n_{\delta}} \sum_{\pm}  \frac{ n |h_{n, \pm}|^2}{R^{-2n} } = 2 \sum_{n = 1}^{n_{\delta}}  R^{2n(1 - \gamma)}
\sim R^{2(1 - \gamma)n_{\delta}} \sim \delta^{-2\alpha}
\end{equation}
and, since $\gamma + 1 > 0$ and $R>1$, 
\begin{equation}\label{tt4}
	\sum_{n = n_{\delta} + 1}^{\infty}  \sum_{\pm} \frac{ n|h_{n, \pm}|^2}{\delta^2 R^{2n}} = \frac{2}{\delta^2} \sum_{n = n_{\delta} + 1}^{\infty}  R^{-2n(\gamma + 1)} \sim  \frac{1}{\delta^2} R^{-2(\gamma + 1)n_\delta} \sim \delta^{-2\alpha}. 
\end{equation}
A combination of \eqref{tt3} and \eqref{tt4} yields \eqref{tt34}. 

It is clear that, since $\gamma > 0$ and $R>1$,  
\begin{equation*}
\| h\|_{H^{1/2}(\partial B_{1/R})}^2 \sim \sum_{n = 1}^\infty \sum_{\pm} n |h_{n, \pm}|^2 = 2 \sum_{n = 1}^\infty R^{-2 n \gamma} < + \infty.
\end{equation*}

\medskip
The proof  is complete. \proofend


\begin{remark} \label{rem-proof-comments} We only prove \eqref{state-3-P2}, \eqref{state-4-P2}, and \eqref{state-5-P2} for the two dimensions. The proof in the three dimensions follows similarly. In fact, in this case, $v_\delta$, $u_\delta$, and $h_\delta$ can be represented by 
\begin{equation*}
\left\{\begin{array}{cl} \dsp v_\delta(x) = \sum_{n = 0}^\infty \sum_{k = -n}^n (a_n^k r^n + b_n^k r^{-n}) Y^k_n(x/ |x|) & \mbox{ in } B_1 \setminus B_{1/R}, \\[6pt]
\dsp u_\delta(x) = \sum_{n = 0}^\infty \sum_{k = -n}^n c_n^k r^n Y^k_n(x/ |x|) & \mbox{ in } B_1, \\[6pt]
\dsp h(x) = \sum_{n = 0}^\infty \sum_{k = -n}^n h_n^k Y^k_n(x/ |x|) &  \mbox{ on } \partial B_{1/R}.  
\end{array}\right.
\end{equation*}
The rest of the proof is almost unchanged.  
\end{remark}

\section{Proof of Theorem~\ref{thm1}} \label{sect-thm1}

\noindent \underline{Step 1:} We show that  if there exists a solution $u \in W^{1}(\mR^{2})$ to 
\begin{equation*}
\dive(\eps_0 \nabla u) = f \mbox{ in } \mR^2, 
\end{equation*}
then $f$ is compatible. This step is not necessary for the proof; however, it gives the motivation for the definition of the  compatibility condition and it guides the proof. 

\medskip 
Define $v$ in $B_{1}$ by
\begin{equation*}
v = u \circ F^{-1}.  
\end{equation*}
We have, by a change of variables, 
\begin{equation}\label{form-changevariables}
 \int_{B_{1} \setminus B_{r}} |\nabla v|^{2}dx
  =  \int_{B_{r^{-1}}\setminus B_{1}} |\nabla u|^{2}dx. 
\end{equation}
Since $v$ is bounded in a neighborhood of the origin, it follows that $v \in H^{1}(B_{1})$ and $\Delta v = F_*f$ in $B_1$ by Lemma~\ref{lem1}.  We have, by Lemma~\ref{lem1} again, 
\begin{equation*}
v = u \Big|_{\rm ext} \mbox{ and } \partial_{r} v = - \partial_{r} u \Big|_{\rm ext} \mbox{ on } \partial B_{1}. 
\end{equation*}
It follows  that 
\begin{equation}\label{transmission-u-P1}
v = u \Big|_{\rm int} \mbox{ and } \partial_{r} v =  \partial_{r} u \Big|_{\rm int} \mbox{ on } \partial B_{1}. 
\end{equation}

Set 
\begin{equation*}
w =  v - u \mbox{ in } B_{1}.
\end{equation*}
Then $w \in H^{1}(B_{1})$ is a solution to the Cauchy problem
\begin{equation*}
	\left\{
		\begin{array}{cl}
			\Delta w =  F_{*}f & \mbox{ in } B_{1},\\[6pt]
			\partial_{\nu} w  = w  = 0 & \mbox{ on } \partial B_1 
		\end{array}
	\right.
\end{equation*}
by \eqref{transmission-u-P1}.  
Therefore $f$ is compatible.

\medskip
\noindent \underline{Step 2:} Proof of statement 1).  

\medskip 
It is clear that $U \in W^{1}(\mR^{2})$ is a solution converging to 0 as $|x| \to \infty$ to  \eqref{eq-P1}.  Statement 1) now follows from the uniqueness of such a solution. 

\medskip
\noindent \underline{Step 3:} Proof of  statement 2).

\medskip 
From \eqref{eq-P1} and \eqref{defTO}, we have 
\begin{equation*}
\int_{B_1} F_* f = 0. 
\end{equation*}
Let $w \in H^1(B_1)$ with $\dsp \int_{B_1} w = 0$ be the unique solution to
\begin{equation}\left\{
	\begin{array}{cl}
		\Delta w = F_{*}f &\mbox{ in }B_1,\\[6pt]
		\partial_{\nu} w = 0 &\mbox{ on }  \partial B_1.
	\end{array}\right.
\end{equation}	 
Define 
\begin{equation}\label{def-Udelta}
	U_{\delta} = \left\{\begin{array}{cl}
		u_{\delta} &\mbox{in } B_1,\\[6pt]
		u_{\delta} - w \circ F &\mbox{in }  \mR^{2} \setminus B_{1}.
	\end{array}\right.
\end{equation}
Similar to  \eqref{form-changevariables}, we have $U_{\delta} \in W^{1}(\mR^{2} \setminus \partial B_1)$. It is clear that
\begin{equation*}
\Delta U_{\delta} = 0 \mbox{ in } \mR^{2} \setminus \partial B_{1}. 
\end{equation*}
Hence, one may represent $U_{\delta}$ as 
\begin{equation}
	U_{\delta} = \left\{\begin{array}{cl}
		\dsp a_0 + \sum_{n = 1}^{\infty} \sum_{\pm} a_{n, \pm} r^n e^{\pm i n \theta} & \mbox{in } B_1,\\[6pt]
		 \dsp b_0 + \sum_{n = 1}^{\infty} \sum_{\pm} b_{n, \pm} r^{-n} e^{\pm i n \theta} & \mbox{in } \mR^{2} \setminus  B_1, 
	\end{array}\right.
\end{equation}
for $a_0, b_0, a_{n, \pm}, b_{n, \pm} \in  \mC$ ($n \ge 1$). Assume that, on $\partial B_{1}$, 
\begin{equation}\label{w-pB1}
w = w_0 + \sum_{n = 1}^{\infty} \sum_{\pm }w_{n, \pm} e^{\pm in \theta}, 
\end{equation}
for $w_0, w_{n, \pm} \in  \mC$ ($n \ge 1$). From \eqref{def-Udelta}, we have
\begin{equation}\left\{
\begin{array}{c}
	U_{\delta} \Big|_{\rm ext} - U_{\delta} \Big|_{\rm int} = -w, \\[6pt]
	\partial_{\nu} U_{\delta} \Big|_{\rm ext} - (-1 + i\delta)\partial_{\nu}U_{\delta} \Big|_{int} = 0.
\end{array}\right.
\end{equation}
This implies 
\begin{equation}
	a_{n, \pm} = b_{n, \pm} + w_{n, \pm} \quad \mbox{ and } \quad (1 - i\delta)a_{n, \pm} = b_{n, \pm}, \quad \, \forall \,  n \geq 1.
\end{equation}
It follows that
\begin{equation}\label{proportional to delta}
	a_{n, \pm} = \frac{w_{n, \pm}}{i\delta} \quad  \mbox{ and } \quad b_{n, \pm} = \frac{(1 - i\delta) w_{n, \pm}}{i\delta} , \quad \, \forall \,  n \geq 1.
\end{equation}
for all $n \geq 1$. Noting that either $w_{n, +} \not = 0$ or $w_{n, -} \not = 0$ for some $n \ge 1$ since  $f$ is not compatible, we obtain \eqref{state-2-P1}. 

\medskip
The proof of Theorem \ref{thm1} is complete. \proofend

\section{Numerical illustrations} \label{Sec Numer}
 
\begin{figure}[h!]
\begin{center}
	\includegraphics[scale=.175]{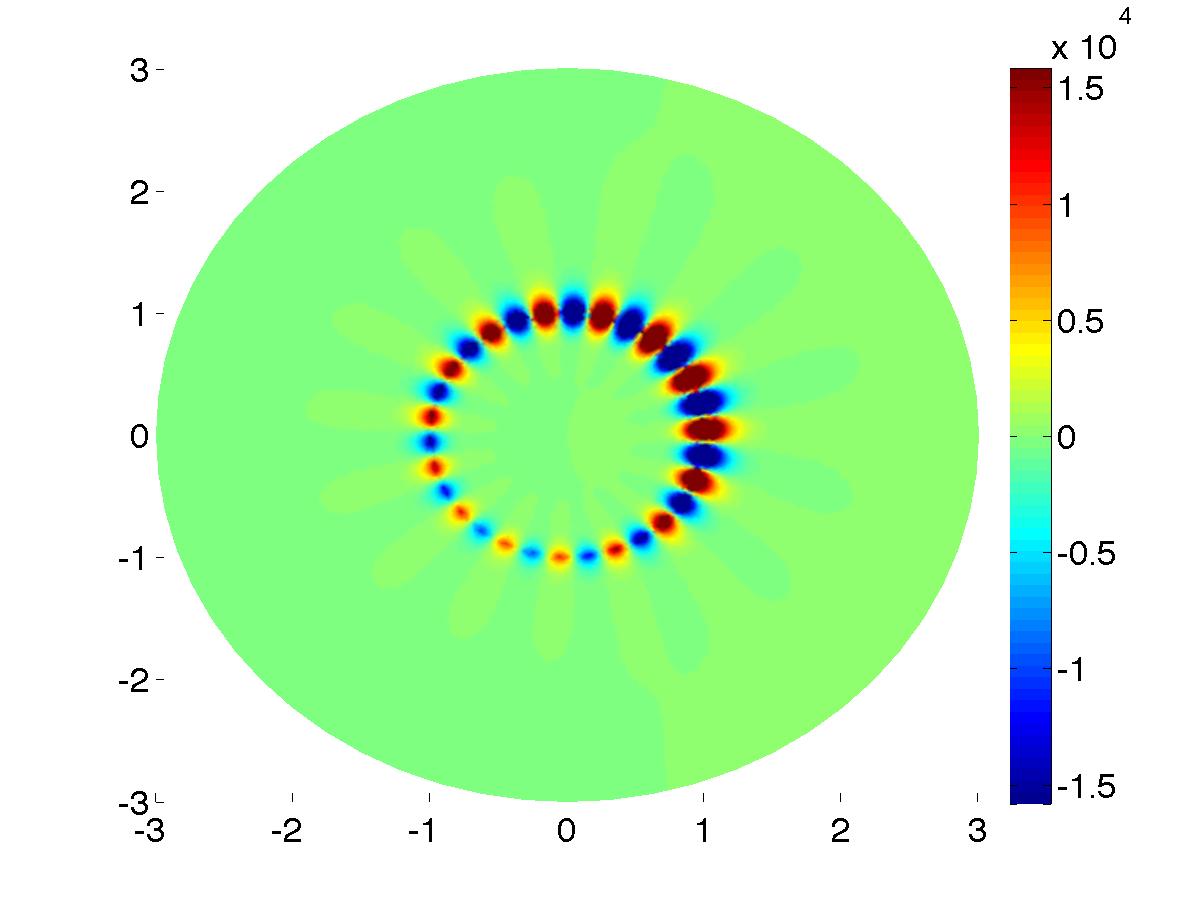} 
	\includegraphics[scale=.175]{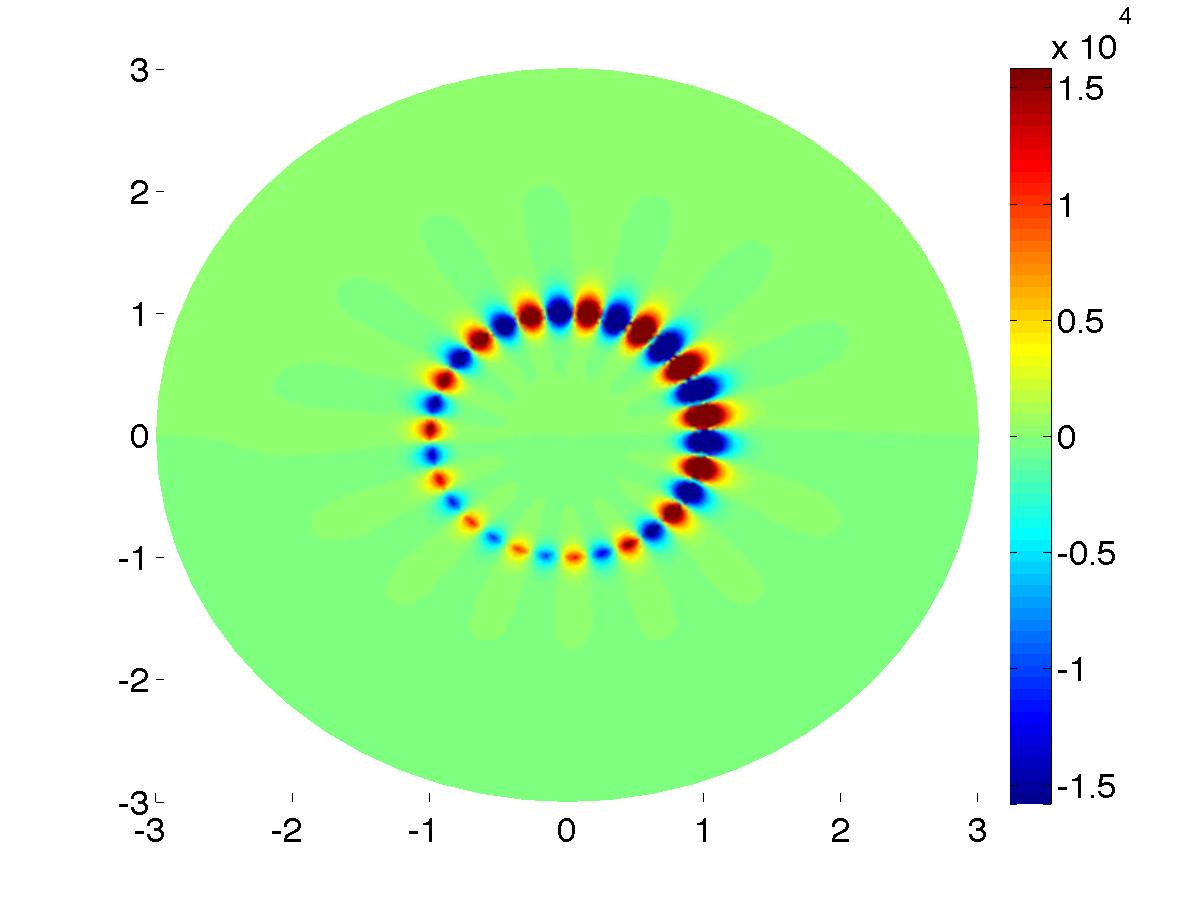} 
	\includegraphics[scale=.175]{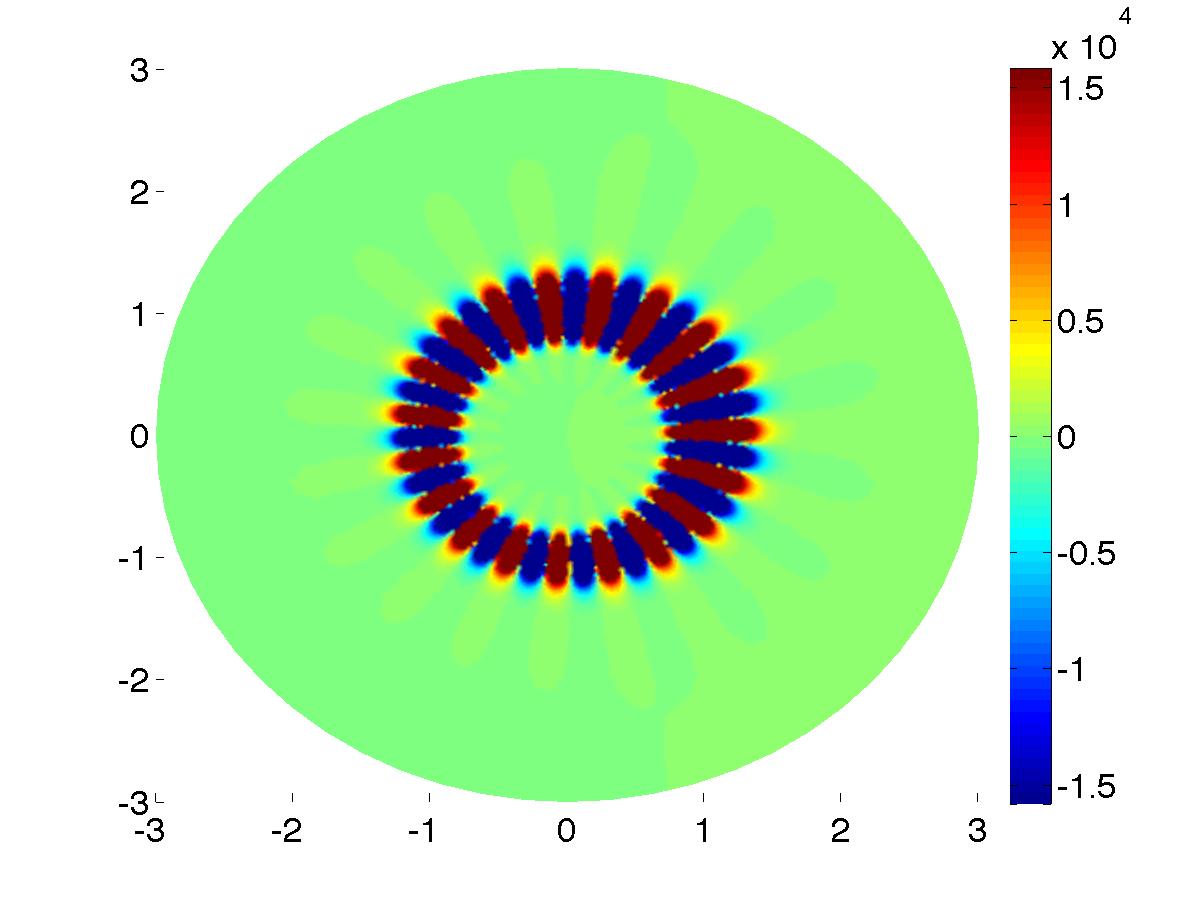} 
	\includegraphics[scale=.175]{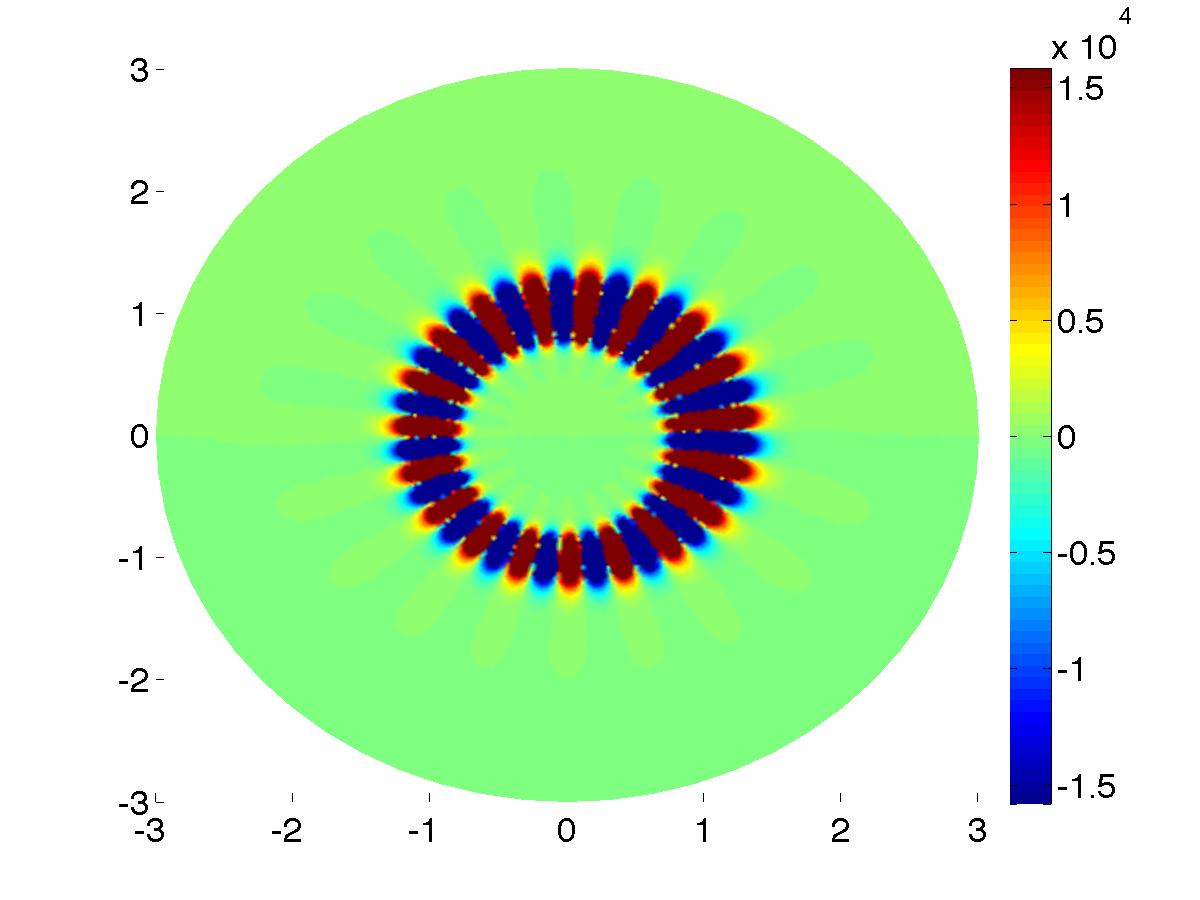}
	\includegraphics[scale=.175]{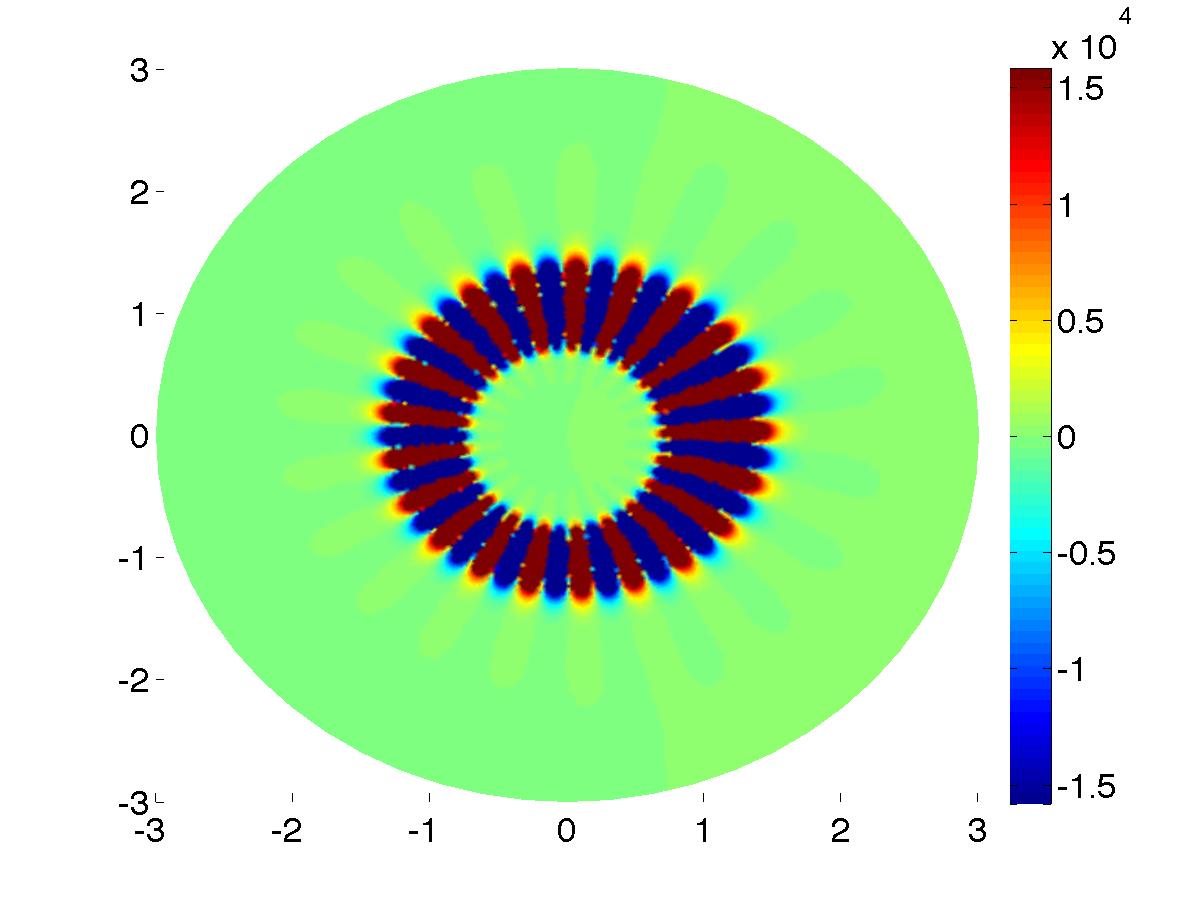} 
	\includegraphics[scale=.175]{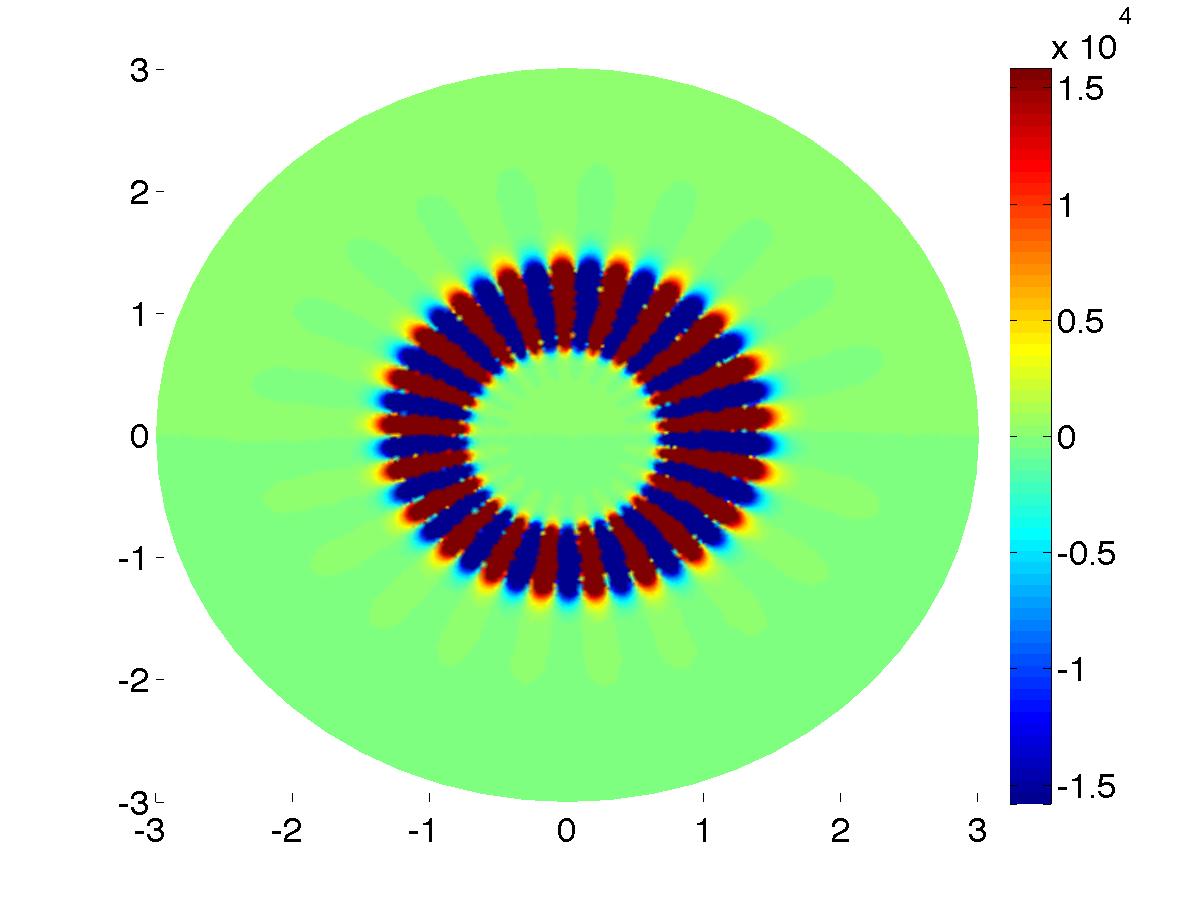} 
\end{center}
\caption{\label{fig 1} The graphs of $u_{\delta}$ when $\delta = 10^{-14}, 10^{-18}$ and $10^{-20}$ from the $1^{\rm st}$ to the $3^{\rm rd}$ row. Left: the real part of $u_{\delta}$; Right: the imaginary part of $u_{\delta}.$}
\end{figure}

\begin{figure}[h!]
\begin{center}
	\includegraphics[scale=.175]{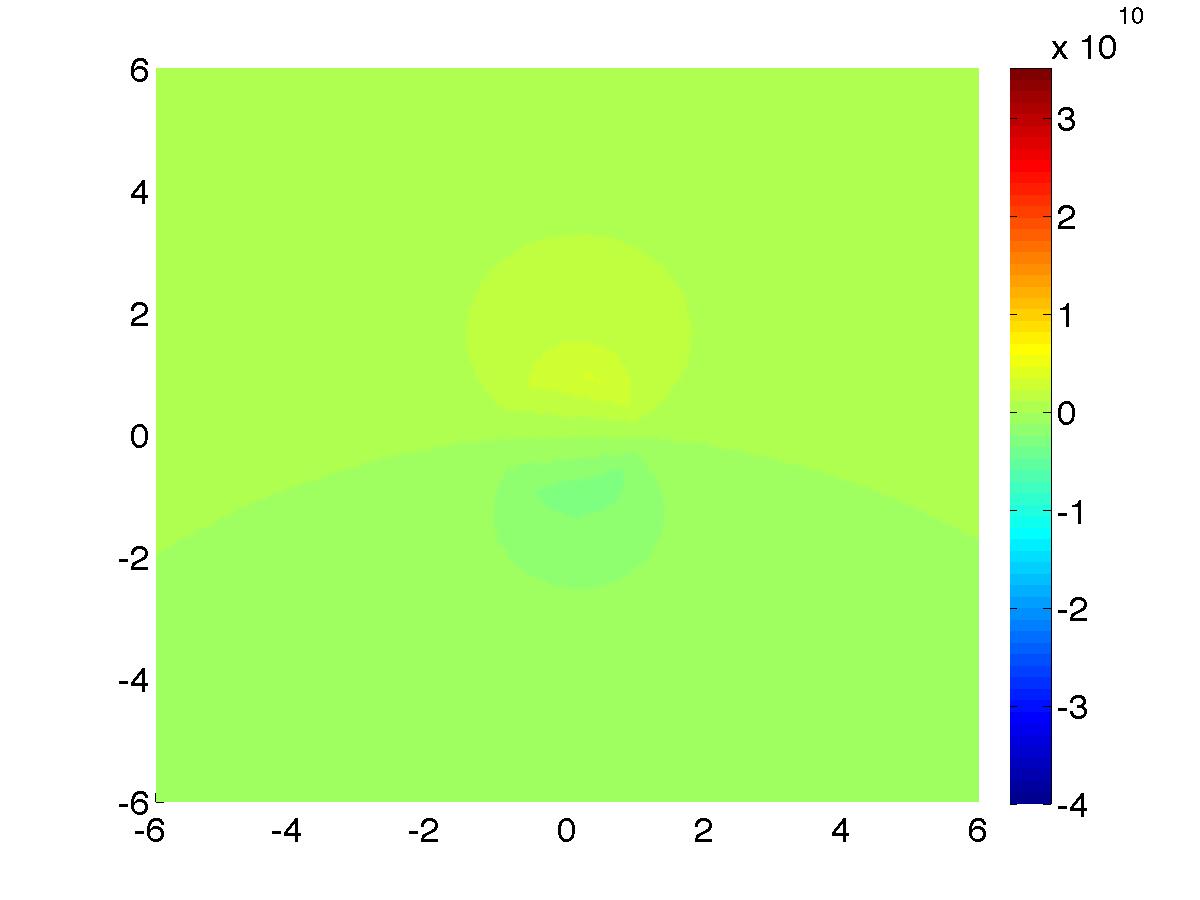}
	\includegraphics[scale=.175]{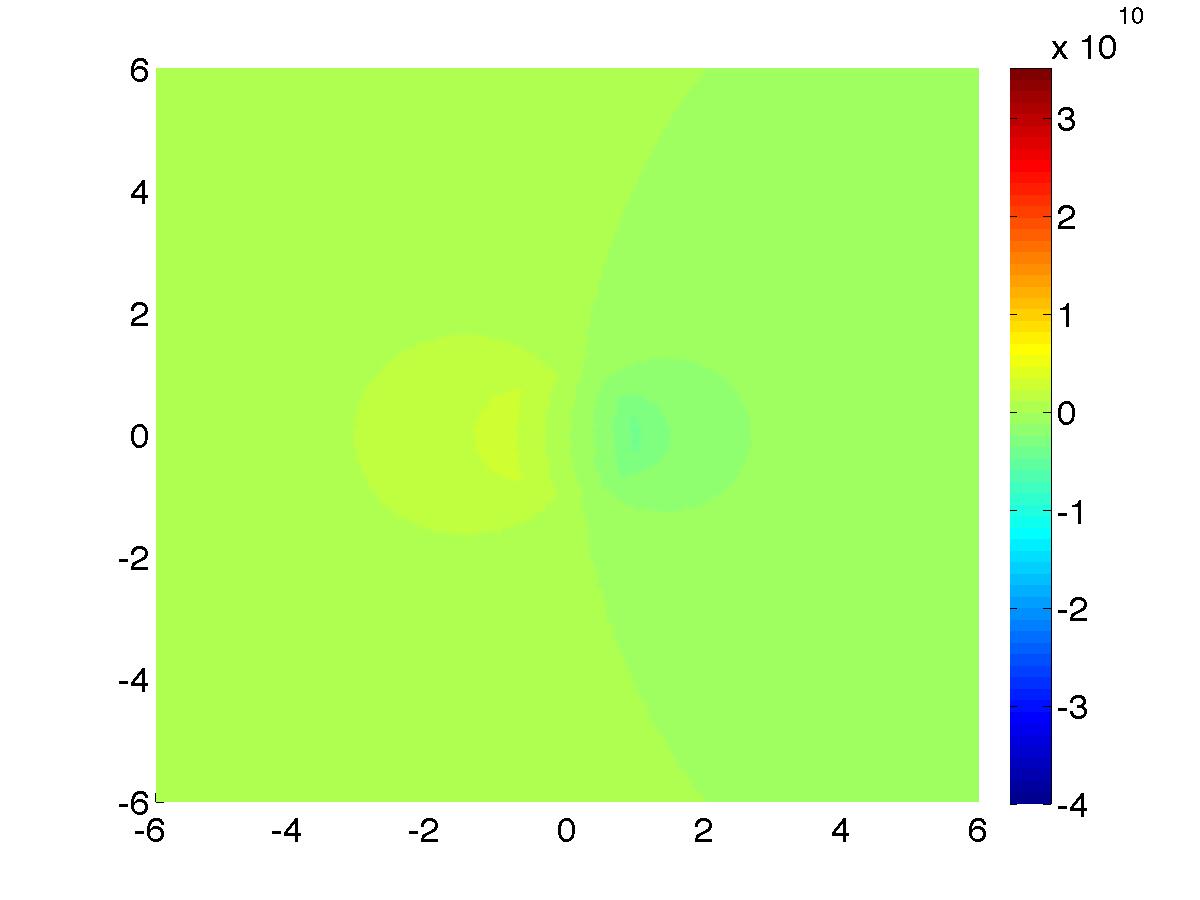} 
	\includegraphics[scale=.175]{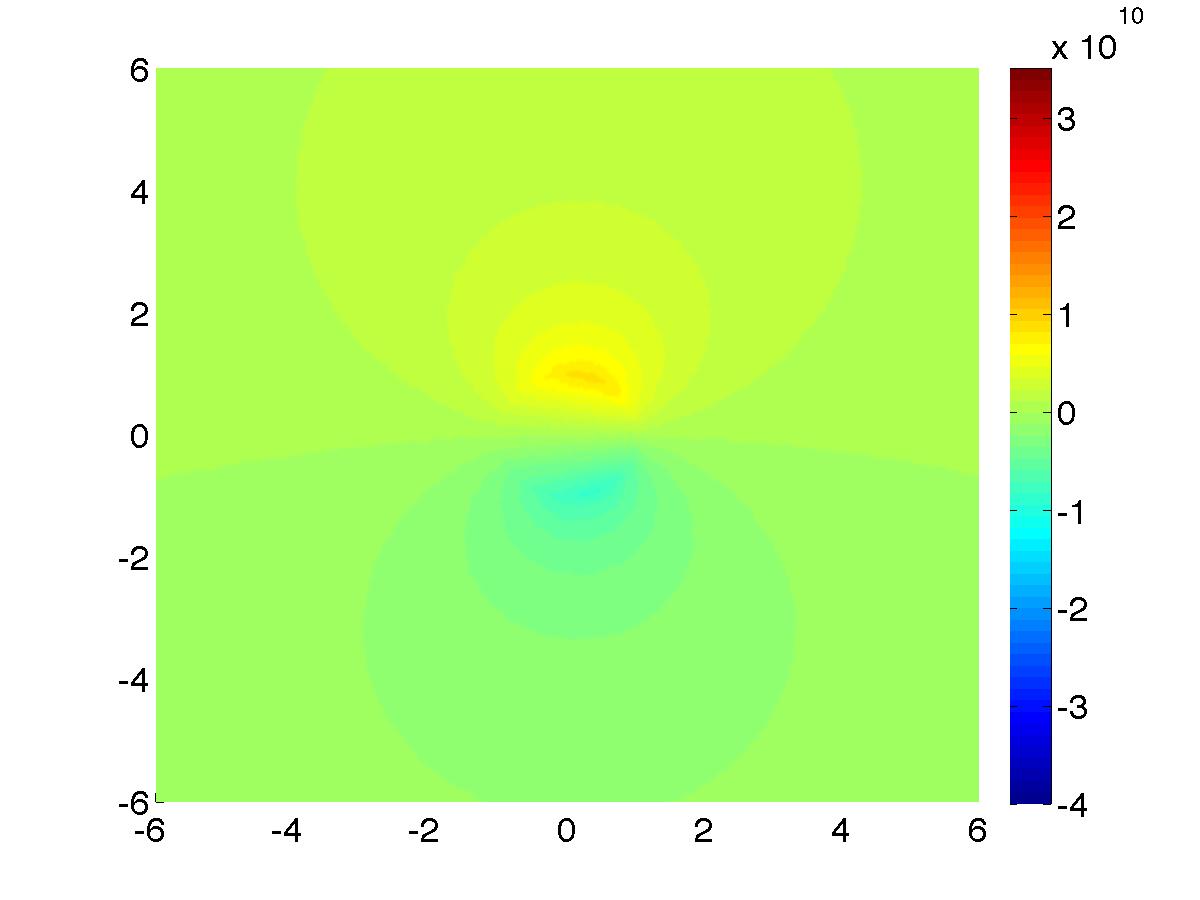} 
	\includegraphics[scale=.175]{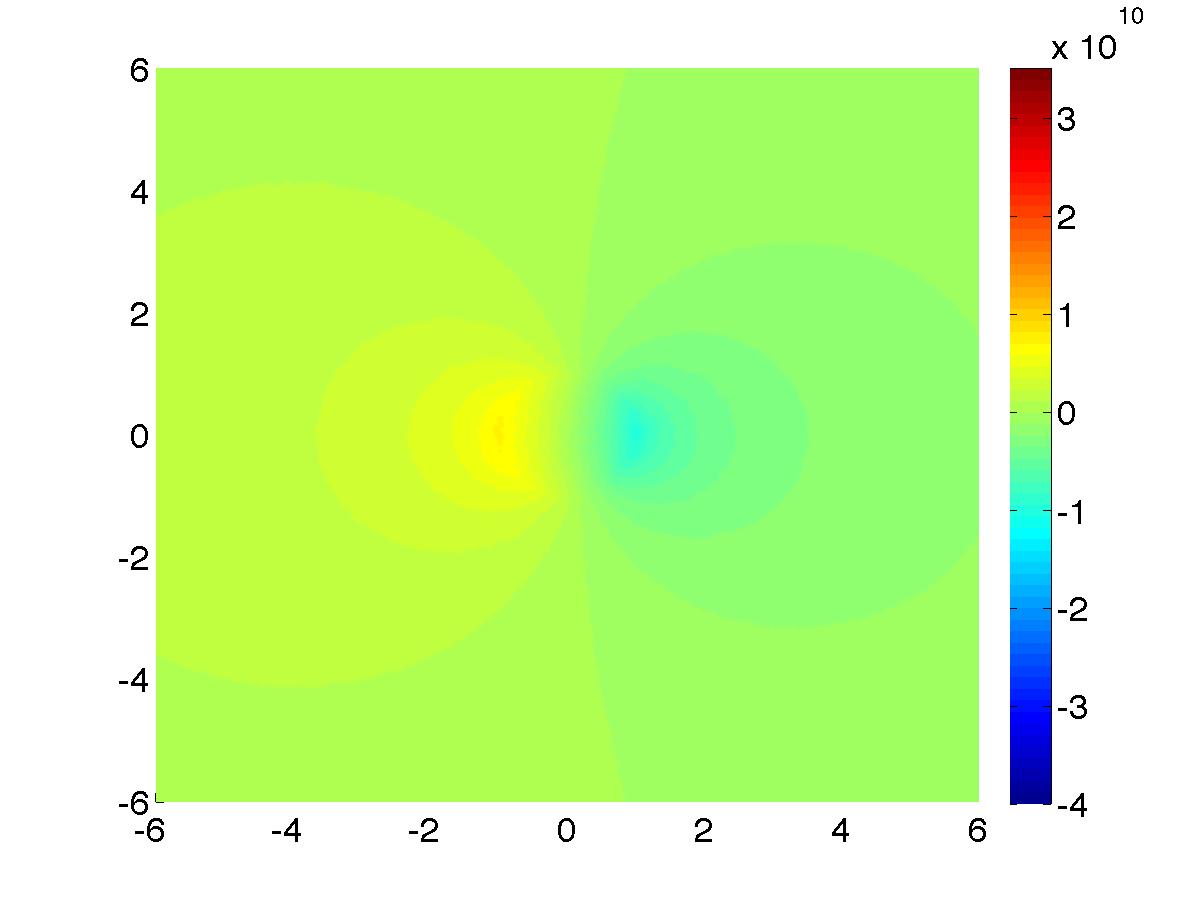} 
	\includegraphics[scale=.175]{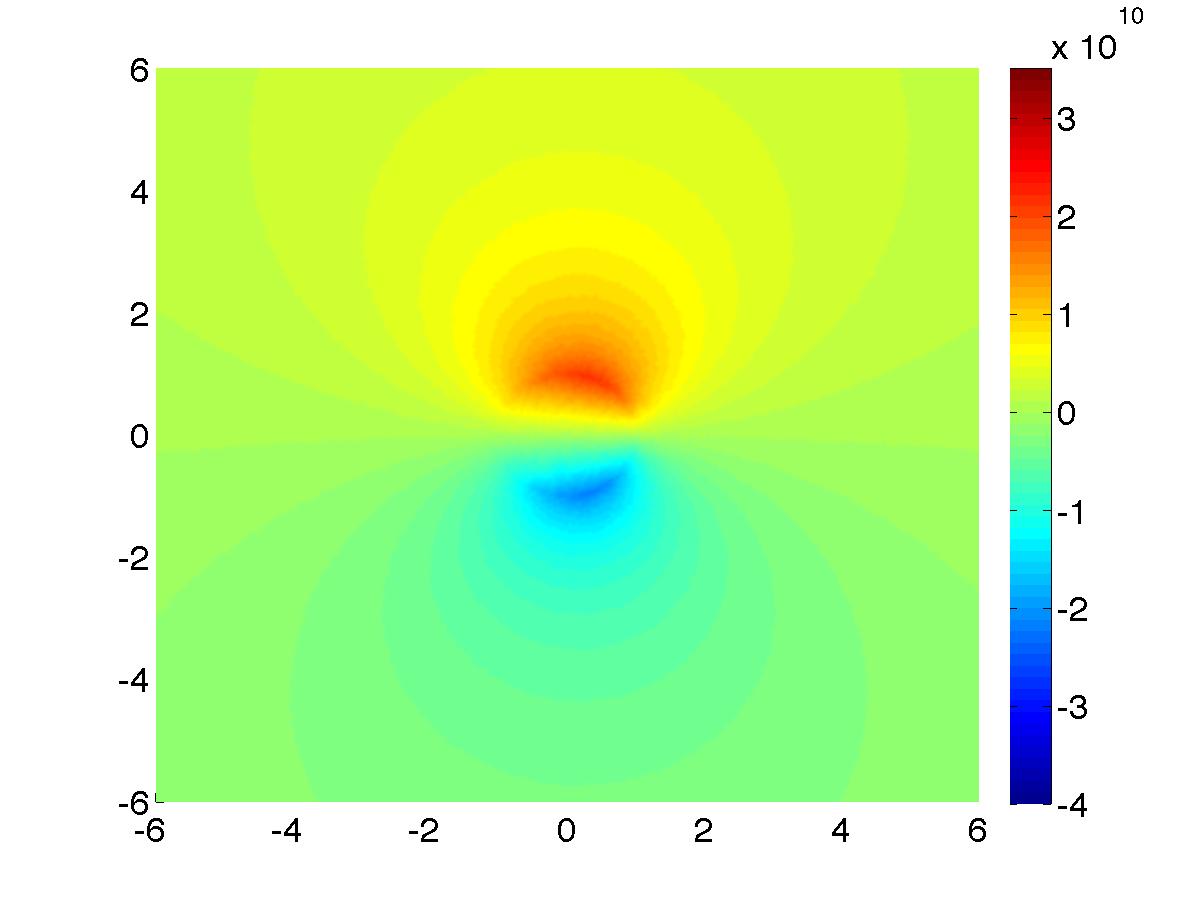}  
	\includegraphics[scale=.175]{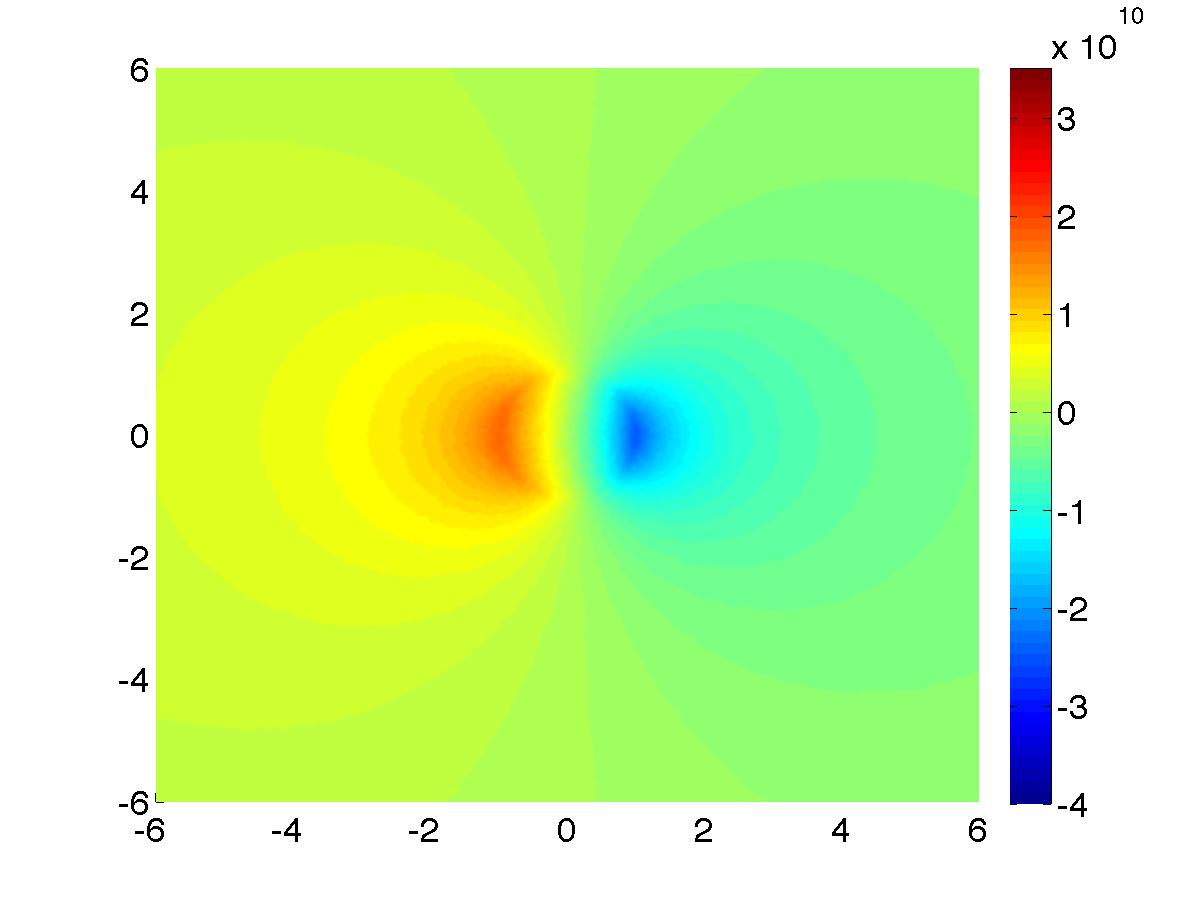} 
\end{center}
\caption{\label{fig 2} The graphs of $u_{\delta}$ when $\delta = 10^{-10}, 10^{-10.4}$ and $10^{-10.8}$ from the $1^{\rm st}$ to the $3^{\rm rd}$ row. Left: the real part of $u_{\delta}$; Right: the imaginary part of $u_{\delta}.$
}
\end{figure}

In this section we present some numerical results to illustrate Theorems \ref{thm2} and \ref{thm1}. Figure \ref{fig 1} corresponds to  Theorem \ref{thm2} and presents a simulation on the localized resonance in which  $R = 3$ and $g =  \sum_{n  = 1}^{\infty} \frac{1}{n^2} e^{i n \theta}$.   Figure \ref{fig 2} corresponds to  Theorem \ref{thm1} and presents a simulation on the complete resonance in which  $f = \Delta (\phi g) \chi_{\mR^2 \setminus B_1}$ where $\chi$ denotes the characteristic function,  $ g = \sum_{n = 1}^{\infty} \frac{r^{n}}{6^n}e^{ i n \theta}$ and $\phi \in C^2(\mR^2)$ is the radially symmetric function such that $\phi = 1$ in $B_2$ and $\phi = 0$ in $\mR^2 \setminus B_3$ \footnote{We take $ \phi(r)  = 513-1080\,r+900\,{r}^{2}-370\,{r}^{3}+75\,{r}^{4}-6\,{r}^{5}$ in $B_3 \setminus B_2$.}. In both simulations, $g$ is approximated by its first hundred terms.

 
\mbox{

}
\newpage
\noindent{\bf Acknowledgement:}  The authors thank the referees for reading 
carefully  the paper and for useful suggestions, in particular  for pointing out the missing point in the proof of Theorems~\ref{thm2} and \ref{thm1} in the previous version.

\providecommand{\bysame}{\leavevmode\hbox to3em{\hrulefill}\thinspace}
\providecommand{\MR}{\relax\ifhmode\unskip\space\fi MR }
\providecommand{\MRhref}[2]{%
  \href{http://www.ams.org/mathscinet-getitem?mr=#1}{#2}
}
\providecommand{\href}[2]{#2}

\end{document}